\documentclass[11pt,a4paper]{article}

\usepackage[latin1]{inputenc}

\usepackage{amsmath}
\usepackage{amsthm}
\usepackage{amsfonts}
\usepackage{latexsym}
\usepackage{amssymb}
\usepackage{amscd}
\usepackage{epsfig}
\usepackage{psfrag}
\usepackage{enumitem}

\usepackage{graphics}
\usepackage{graphicx}
\usepackage{bezier}
\usepackage{float}
\usepackage{color}
\usepackage{hyperref}
\usepackage{url}

\usepackage{multienum}
\usepackage[table]{xcolor}
\usepackage{multicol}
\usepackage{multirow}
\usepackage{fancyvrb}
\usepackage{parskip}
\usepackage[toc,page]{appendix}
\usepackage{blkarray}
\usepackage[none]{hyphenat}

\usepackage[
top=2.7cm,
bottom=2.7cm,
left=1.5cm,
right=1.5cm
]{geometry}

\newtheorem{theorem}{Theorem}[section]
\newtheorem{lemma}{Lemma}[section]
\newtheorem{proposition}{{Proposition}}[section]

\newtheorem{remark}{Remark}[section]

\numberwithin{equation}{section}
\numberwithin{table}{section}
\numberwithin{figure}{section}

\title{On Wagstaff primes in the $k$-Pell number sequence}
\author{Prosper Kaggwa $^{1,*}$}
\date{}

\begin{document}
	
	\maketitle
\begin{center}
	\small
	\textit{In memory of Professor Florian Luca (1969--2026)}
\end{center}
	\begin{abstract}
		A Wagstaff prime is a prime number of the form $(2^{\mathfrak{p}}+1)/3$, where $\mathfrak{p}$ is an odd prime.
		For an integer $k\ge 2$, the sequence of $k$-generalized Pell numbers $\bigl(P_n^{(k)}\bigr)_{n\ge 2-k}$
		is the linear recurrence of order $k$ given by
		$P_n^{(k)}=2P_{n-1}^{(k)}+P_{n-2}^{(k)}+\cdots+P_{n-k}^{(k)}$ for $n\ge 2$,
		with initial values $P_{2-k}^{(k)}=\cdots=P_{0}^{(k)}=0$ and $P_1^{(k)}=1$.
		For $k=2$ this is the classical Pell sequence.
		Motivated by the recent determination of the Wagstaff primes in the $k$-Fibonacci and $k$-Lucas sequences,
		we study the exponential Diophantine equation
		$P_n^{(k)}=(2^{\mathfrak{p}}+1)/3$ in positive integers $n,k$ and an odd prime $\mathfrak{p}\geq 3$, with $k\ge 2$.
		We show that this equation has no solutions, that is no $k$-Pell number is a Wagstaff prime.
		The proof combines linear forms in logarithms of algebraic numbers
		with reduction based on the theory of continued fractions, together with a computational verification in the remaining finite range.
	\end{abstract}
	
	
	
	\thanks{$ ^{*} $ Corresponding author}
	
	\section{Introduction}\label{intro}
	\subsection{Background}
	The search for terms of a linear recurrence sequence that possess a prescribed arithmetic shape is a classical theme of Diophantine Number theory. Over the last two decades, the tools of Baker's theory of linear forms in logarithms, in the explicit form due to \cite{Matveev2000} (see also \cite{Bugeaud2006}), combined with reduction methods such as the Baker--Davenport lemma in the version of \cite{DP} have made it possible to determine completely the members of families of sequences such as generalized Fibonacci and Lucas sequences that are, for instance, powers of two \cite{Bravo_Luca2012}, repdigits
	\cite{Bravo_Luca2014}, Mersenne numbers \cite{Bravo_Gomez2016,Altassan_Alan2024_Mersenne}, or Fermat numbers \cite{Bravo_Herrera2020}.
	
	A \emph{Wagstaff prime} is a prime number of the form
	\[
	W_{\mathfrak{p}}=\frac{2^{\mathfrak{p}}+1}{3},
	\]
	where $\mathfrak{p}$ is an odd prime. The first few Wagstaff primes are $3, 11, 43, 683, 2731, 43691,\dots$, corresponding to
	$\mathfrak{p}=3,5,7,11,13,17,\dots$ (see sequence A000979 in \cite{Sloane_A000979}). They are named after Samuel S.\ Wagstaff Jr.\ and
	are closely tied to Mersenne primes through the New Mersenne Conjecture of \cite{Bateman_Selfridge_Wagstaff1989}. It is not known whether
	there are infinitely many Wagstaff primes, and their primality has been studied both through compositeness tests
	\cite{Berrizbeitia_Luca_Melham2010} and through the search for probable primes \cite{Sikkema2024}.
	
	Very recently, Wagstaff primes have been sought inside generalized Fibonacci-type sequences. For $k\ge 2$, let
	$\bigl(F_n^{(k)}\bigr)_{n\ge 2-k}$ be the $k$-generalized Fibonacci sequence, given by
	$F_n^{(k)}=F_{n-1}^{(k)}+\cdots+F_{n-k}^{(k)}$ with initial values
	$F_{2-k}^{(k)}=\cdots=F_0^{(k)}=0$ and $F_1^{(k)}=1$.
	\cite{Rezaiguia_Qawaqneh_Abdelouahab2026} solved the equation $F_n^{(k)}=(2^{\mathfrak{p}}+1)/3$ and proved that the only Wagstaff prime
	appearing in any $k$-Fibonacci sequence is $3$, namely $F_4^{(2)}=3$. Shortly afterwards, \cite{Batte2026}
	treated the $k$-Lucas sequence $\bigl(L_n^{(k)}\bigr)_{n\ge 2-k}$, whose initial values are $L_0^{(k)}=2$, $L_1^{(k)}=1$
	and $L_j^{(k)}=0$ for $2-k\le j\le -1$, and showed that the equation $L_n^{(k)}=(2^{\mathfrak{p}}+1)/3$ has only the solutions
	$(n,k,\mathfrak{p})\in\{(5,2,5),(6,4,7)\}\cup\{(2,k,3):k\ge 2\}$; in particular the Wagstaff primes $3$, $11$ and $43$
	occur in this family.
	
	The present paper considers the same question for the \emph{$k$-generalized Pell numbers}, studied by
	\cite{Bravo_Herrera_Luca2021}. For an integer $k\ge 2$, this is the sequence $\bigl(P_n^{(k)}\bigr)_{n\ge 2-k}$ defined by the recurrence
	\[
	P_n^{(k)}=2P_{n-1}^{(k)}+P_{n-2}^{(k)}+\cdots+P_{n-k}^{(k)},\qquad n\ge 2,
	\]
	with initial conditions $P_{2-k}^{(k)}=P_{3-k}^{(k)}=\cdots=P_{0}^{(k)}=0$ and $P_1^{(k)}=1$.
	For $k=2$ we recover the Pell numbers. Although the $k$-Pell sequences resemble the $k$-Fibonacci and $k$-Lucas sequences, the argument cannot be transferred verbatim from the previous two papers.
	
	Let $\mathfrak{p}$ be an odd prime. The main result of this paper is the following.
	\subsection{Main Result} 
	\begin{theorem}\label{thm:main}
		Let $\bigl(P_n^{(k)}\bigr)_{n\ge 2-k}$ be the sequence of $k$-generalized Pell numbers. The Diophantine equation
		\begin{equation}\label{eq:main}
			P_n^{(k)}=\frac{2^{\mathfrak{p}}+1}{3}
		\end{equation}
		has no solutions in integers $n\ge 1$, $k\ge 2$ and primes $\mathfrak{p}\geq 3$. 
	\end{theorem}
	\section{Methods}
	\subsection{Preliminaries}
	The following identity, taken from \cite{kilis}, relates $P_n^{(k)}$ to the $m$-th Fibonacci number where $m=2n-1$. It states that
	
	\begin{align}\label{eq:2.1}
		P_n^{(k)} =  F_{2n-1},
	\end{align}
	
	for all $1 \le n \le k+1$, while the next term satisfies $P_{k+2}^{(k)} =  F_{2k+3}-1$. The characteristic polynomial of the sequence $(P_n^{(k)})_{n\geq 2-k}$ is
	\[
	\Psi_k(x) = x^k - 2x^{k-1} - \cdots - x - 1.
	\]
	
	This polynomial is irreducible over $\mathbb{Q}[x]$ and has a unique real root $\alpha := \alpha(k) > 1$, with all other roots lying inside the unit circle, see \cite{Bravo_Herrera_Luca2021}. Moreover, the root $\alpha$ satisfies
	
	\begin{align}\label{eq2.3m}
		\varphi^2(1 - \varphi^{-k}) < \alpha < \varphi^2, \qquad k \ge 2, \quad \text{where } \varphi := \frac{1+\sqrt{5}}{2},
	\end{align}
	
	as noted in \cite{Bravo_Herrera_Luca2021}. Similarly to the case $k=2$, it has been shown in \cite{Bravo_Herrera_Luca2021} that
	
	\begin{align}\label{fn_up}
		\alpha^{n-2} \le P_n^{(k)} \le \alpha^{n-1}, \qquad \text{for all } n \ge 1, \, k \ge 2.
	\end{align} For $k \ge 2$, we define
	\begin{align*}
		f_k(x) := \frac{x - 1}{(k + 1)x^2 - 3kx + k - 1} = \frac{x - 1}{k(x^2 - 3x + 1) + x^2 - 1}.
	\end{align*}
	
	Bravo and Herrera \cite{BraHer} showed in Lemma 1 that
	\begin{align*}
		0.276 < f_k(\alpha) < 0.5 \qquad \text{and} \qquad |f_k(\alpha_i)| < 1,
	\end{align*}
	for all $2 \le i \le k$, where $\alpha_i$ are the roots of $\Psi_k(x)$ inside the unit circle. This shows that $f_k(\alpha)$ is not an algebraic integer. Furthermore, Theorem 3.1 in \cite{Bravo_Herrera_Luca2021} establishes
	\begin{align}\label{eq2.6m}
		P_n^{(k)} = \sum_{i=1}^{k} f_k(\alpha_i)\alpha_i^n \qquad \text{and} \qquad \left| P_n^{(k)} - f_k(\alpha)\alpha^n \right| < \frac{1}{2},
	\end{align}
	for all $k \ge 2$ and $n \ge 2-k$. The first part of \eqref{eq2.6m} provides a Binet-type formula for $P_n^{(k)}$, while the second inequality shows that the influence of the roots inside the unit circle is negligible.
	
	A sharper estimate appears in Lemma 2 of \cite{BraHer}, valid under a more restricted range of $n$ relative to $k$. Specifically, if $k \ge 30$ and $n > 1$, then
	\begin{align}\label{pk_b1}
		\left| f_k(\alpha)\alpha^n - \frac{\varphi^{2n}}{\varphi+2} \right| < \left( \frac{\varphi^{2n}}{\varphi+2} \right) \frac{4}{\varphi^{k/2}}, \qquad \text{provided } n < \varphi^{k/2}.
	\end{align}
	
	\subsection{Linear Forms in Logarithms}
	
	Let \(\gamma\) be an algebraic number of degree \(d\), and let its minimal polynomial over \(\mathbb{Z}\) be  
	\[
	u_0 x^d + u_1 x^{d-1} + \cdots + u_d = u_0 \prod_{i=1}^{d} \bigl(x - \gamma^{(i)}\bigr),
	\]  
	where the coefficients \(u_i\) are coprime integers with \(u_0 > 0\), and \(\gamma^{(i)}\) denote the conjugates of \(\gamma\).  
	The logarithmic height of \(\gamma\) is defined as  
	\[
	h(\gamma) = \frac{1}{d}\Biggl( \log u_0 + \sum_{i=1}^{d} \log\max\bigl\{|\gamma^{(i)}|, 1\bigr\} \Biggr).
	\]
	
	In particular, for a rational number \(r/s\) with \(\gcd(r,s) = 1\) and \(s > 0\), the logarithmic height is \(h(r/s) = \log\max\{|r|, s\}\).  
	The following standard properties of logarithmic heights will be applied without further reference:
	\begin{enumerate}[label=(\alph*)]
		\item \(h(\gamma_1 \pm \gamma_2) \le h(\gamma_1) + h(\gamma_2) + \log 2\),
		\item \(h(\gamma_1\gamma_2^{\pm 1}) \le h(\gamma_1) + h(\gamma_2)\),
		\item \(h(\gamma^{s}) = |s|\, h(\gamma)\) for any integer \(s\).
	\end{enumerate}
	
	Based on these properties, it follows from \cite[Lemma 3]{Bravo_Herrera_Luca2021} that
	\begin{align}\label{eqh}
		h\bigl(f_k(\alpha)\bigr) < 4k\log\varphi + k\log(k+1), \qquad \text{for all } k \ge 2.
	\end{align}
	
	The following theorem, proved in \cite{Matveev2000} and cited in \cite[Theorem 9.4]{Bugeaud2006}, is fundamental in bounding variables in the Diophantine equation \eqref{eq:main}.  
	Let \(\mathbb{K}\) be a number field of degree \(D\), and let \(\gamma_1, \dots, \gamma_t \in \mathbb{K}\) be non-zero elements with integers \(b_1, \dots, b_t\).  
	Define
	\[
	B = \max\{|b_1|, \dots, |b_t|\}, \qquad 
	\Gamma = \gamma_1^{b_1} \cdots \gamma_t^{b_t} - 1,
	\]
	and let \(A_1, \dots, A_t\) be real numbers satisfying
	\[
	A_j \ge h'(\gamma_j) := \max\{D\,h(\gamma_j),\; |\log\gamma_j|,\; 0.16\}, \quad 1 \le j \le t.
	\]
	We call \(h'\) the modified height (with respect to the field \(\mathbb{K}\)).  
	With this notation, the main result in \cite{Matveev2000} implies the following estimate.
	
	\begin{theorem}[Theorem 9.4 in \cite{Bugeaud2006}]\label{matveev}
		Assume that \(\Gamma\) is non-zero and \(\mathbb{K}\) is real. Then
		\[
		\log|\Gamma| > -1.4 \cdot 30^{\,t+3}\, t^{4.5}\, D^{2}\, 
		(1+\log D)(1+\log B)\, A_1\cdots A_t.
		\]    
	\end{theorem}
	The computations initially lead to upper bounds for the variables that are too large to be used directly. We therefore apply a Baker-Davenport type reduction to obtain smaller bounds. The details of this reduction procedure are presented in the following subsection. The precise variant we need is due to Dujella and Peth\H{o} \cite{DP}.
	\subsection{Reduction method}\label{sec2.3}  
	For a real number \(r\) we write \(\|r\|\) for the distance from \(r\) to the nearest integer, i.e.
	\(\|r\|:=\min\{|r-n|:n\in\mathbb{Z}\}\).
	
	\begin{lemma}[Dujella, Peth\H{o}, \cite{DP}]\label{DP}
		Let \(\tau\neq0\), \(C\), \(\beta\) and \(\mu\) be real numbers with \(C>0\) and \(\beta>1\).
		Let \(M>1\) be a positive integer and suppose that \(r/s\) is a convergent of the continued fraction
		expansion of \(\tau\) satisfying \(s>6M\).  Put
		\[
		\theta:=\|\mu s\|-M\|\tau s\|.
		\]
		If \(\theta>0\), then the inequality
		\[
		0<|u\tau-v+\mu|<C\beta^{-\ell}
		\]
		has no solution in positive integers \(u,v,\ell\) with
		\[
		\frac{\log(Cs/\theta)}{\log \beta}\le \ell \qquad\text{and}\qquad u\le M .
		\]
	\end{lemma}
	The following is an important fact from calculus that will be useful in the proof.\\
	If $z \in \mathbb{R}$ satisfies $|z|<\frac{1}{2}$, then 
	\begin{align}\label{eq.2,7}
		\notag	|\log{(1+z)}|
			&=\left|z-\frac{z^2}{2}+\frac{z^3}{3}-\frac{z^4}{4}+\cdots \right|,\\
			&<|z|+\frac{1}{2}\left(|z|^2+|z|^3+|z|^4+\cdots \right)
			=|z|\left(1+\frac{|z|}{2(1-|z|)}\right)
			<\frac{3}{2}|z|.
	\end{align}
	Finally, we present an analytic argument which is Lemma 7 in \cite{GL}.
	\begin{lemma}[Lemma 7 in \cite{GL}]\label{Guz}
		If $r \geq 1$, $T > (4r^2)^r$ and $T > \dfrac{x}{(\log x)^r}$, then
		\[
		x < 2^r T (\log T)^r.
		\]
	\end{lemma}	SageMath 10.6 is used to perform all computations in this work.
	\section{Proof of Theorem 1.1}
	\subsection{Some necessary bounds}
	For every odd prime $\mathfrak{p}$, equation \eqref{eq:main} together with the lower
	estimate for $P_n^{(k)}$ in terms of $\alpha$ gives
	\[
	\alpha^{n-2} < P_n^{(k)} = \frac{2^{\mathfrak{p}}+1}{3} < 2^{\mathfrak{p}}+1 < 2^{\mathfrak{p}+1}.
	\]
	Taking logarithms, we get $(n-2)\log\alpha < (\mathfrak{p}+1)\log 2$, so that
	$n < 2\mathfrak{p}+4$. 
	In the other direction, the upper estimate for $P_n^{(k)}$ gives
	\[
	\alpha^{n-1} \ge P_n^{(k)} = \frac{2^{\mathfrak{p}}+1}{3} > \frac{2^{\mathfrak{p}}}{3} > 2^{\mathfrak{p}-2},
	\]
	and taking logarithms once more yields $(n-1)\log\alpha \ge (\mathfrak{p}-2)\log 2$, since $\log {\alpha}<2\log 2$ from $\alpha<\varphi^2<4,$ we get $n > \mathfrak{p}/2$. Putting the two estimates together,
	\begin{equation}\label{eq:n-range}
		\frac{\mathfrak{p}}{2} < n < 2\mathfrak{p}+4.
	\end{equation}
	
	\subsection{The case $n \le k+1$}
	\begin{proposition}
		Let $k\geq 2, 1\leq n\leq k+1$ and let $\mathfrak{p}\geq 3$ be an integer. Then $P_n^{(k)}\neq \left(2^{\mathfrak{p}}+1\right)/3.$
	\end{proposition}
	\begin{proof}
		By \eqref{eq:2.1}, $P_n^{(k)} = F_{2n-1}$. Put $W = (2^{\mathfrak{p}}+1)/3$. Since $\mathfrak{p} \ge 3$,
		$8 \mid 2^{\mathfrak{p}}$, so $3W = 2^\mathfrak{p} + 1 \equiv 1 \pmod{8}$. As
		$3 \cdot 3 \equiv 1 \pmod{8}$, we get $W \equiv 3 \pmod{8}$. The Fibonacci
		sequence modulo $8$ is periodic with period $12$; one period is
		$0,1,1,2,3,5,0,5,5,2,7,1$ (indices $m = 0,\ldots,11$), so
		$F_m \equiv 3 \pmod{8}$ only if $m \equiv 4 \pmod{12}$, in particular only
		if $m$ is even. Since $2n-1$ is odd, $F_{2n-1} \not\equiv 3 \pmod{8}$, and
		therefore $F_{2n-1} \ne W$.
	\end{proof}
	\begin{remark}
		For $k=2$ the Pell numbers satisfy $P_n^{(2)} \bmod 8 \in \{0,1,2,5,4,5,6,1\}$
		(period $8$) and are never $\equiv 3 \pmod{8}$; so Theorem~\ref{thm:main} is already
		complete for $k=2$ by the same argument, without any assumption on $n$.
		For $k \ge 3$ the residue $3$ does occur, so the analytic argument of Subsection~3.3 is needed.
	\end{remark}
	
	We therefore have that equation \eqref{eq:main} has no solution with
	$n \le k+1$ and $k \ge 2$, and from now on we assume $n \ge k+2$.
	
	\subsection{The case $n \ge k+2$}
	
	Here $n \ge k+2$, and since $k \ge 2$ this forces $n \ge 4$.
	We first bound $n$ in terms of $k$.
	
	\begin{lemma}\label{lem:n-bound}
		Let $(n,k,\mathfrak{p})$ be an integer solution of \eqref{eq:main} with $n \ge k+2$,
		$k \ge 2$ and $\mathfrak{p} \ge 3$ prime, then
		\begin{equation}\label{eq:n-bound}
			n < 1.3\cdot 10^{15}\, k^{6} (\log k)^{3}.
		\end{equation}
	\end{lemma}
	
	\begin{proof}
	We return to \eqref{eq:main} and use the closed formula \eqref{eq2.6m}. Write
	$$P_n^{(k)}-f_k(\alpha)\alpha^{n}=e_k(n),$$ such that $|e_k(n)|<1/2$. Then
	\[
	\frac{2^{\mathfrak{p}}+1}{3}-f_k(\alpha)\alpha^{n}=e_k(n),\] thus
	\[\frac{2^{\mathfrak{p}}}{3}-f_k(\alpha)\alpha^{n}=e_k(n)-\frac{1}{3}.\]
	Taking absolute values on both sides yields
	\[
	\left|\frac{2^{\mathfrak{p}}}{3}-f_k(\alpha)\alpha^{n}\right|
	\le |e_k(n)|+\frac13<\frac12+\frac13=\frac56 .
	\]	Thus we have
		\[
		\left| \frac{2^{\mathfrak{p}}}{3} - f_k(\alpha)\alpha^{n} \right|
		< \frac{5}{6}.\]
		Dividing through by $f_k(\alpha) \alpha^n$ yields
		\[
		\left| 2^{\mathfrak{p}} \bigl( 3 f_k(\alpha)\alpha^{n} \bigr)^{-1} - 1 \right|
		<\frac{5}{6f_k(\alpha)\alpha^n}< \frac{4}{\alpha^{n}}.
		\] where we have used the fact that $f_k(\alpha)>0.276$.
		
		Set
		\begin{equation}\label{eq:Gamma-bound}
			|\Gamma_1| := \left| 2^{\mathfrak{p}} \bigl( 3 f_k(\alpha)\alpha^{n} \bigr)^{-1} - 1 \right|
			< \frac{4}{\alpha^{n}}.
		\end{equation}
		Clearly, $\Gamma_1 \ne 0$ otherwise we would have
		$2^{\mathfrak{p}} = 3 f_k(\alpha)\alpha^{n}$. Applying a Galois automorphism $\sigma$ of the Galois closure of
		$\mathbb{Q}(\alpha)$, with $\sigma(\alpha)=\alpha_i$ for some
		$i>1$, and taking absolute values on both sides, we obtain
		\[
		2^{\mathfrak{p}} = 3\,\bigl| f_k(\alpha_i)\alpha_i^{n} \bigr|.
		\]
		This is a contradiction, since the left-hand side is at least $8$ while the
		right-hand side is less than $3$. 
		\begin{remark}
			The above argument is valid for every $\mathfrak{p}\geq 3$.
		\end{remark}
		
		To apply Theorem~\ref{matveev}, we work in the field $\mathbb{K} := \mathbb{Q}(\alpha)$,
		so $D = k$, and we take $t = 3$ with
		\begin{align*}
			\gamma_1 &:= 2, & \gamma_2 &:= 3 f_k(\alpha), & \gamma_3 &:= \alpha,\\
			b_1 &:= \mathfrak{p}, & b_2 &:= -1, & b_3 &:= -n.
		\end{align*}
		We have $h(\gamma_1) = h(2) = \log 2$, 
		\[
		h(\gamma_3) = h(\alpha) < \frac{\log\alpha}{k} < \frac{2\log\varphi}{k} < \frac{1}{k},
		\]
		and for $\gamma_2$,
		\begin{align*}
			h(\gamma_2) &= h\bigl(3 f_k(\alpha)\bigr),\\
			&\leq h(3) + h\bigl(f_k(\alpha)\bigr),\\
			&< \log 3 + 4k\log\varphi + k\log(k+1),\\
			&< 6k\log k.
		\end{align*}
		We therefore take $A_1 := k\log 2$, $A_2 := 6k^{2}\log k$ and
		$A_3 := k$.
		It remains to choose $B \ge \max\{|b_i| : i = 1,2,3\}$; by
		\eqref{eq:n-range} we take $B := 2n$.
		
		Theorem~\ref{matveev} now gives
		\begin{align}
			\log|\Gamma_1| &> -1.4\cdot 30^{6}\cdot 3^{4.5}\cdot k^{2}(1+\log k)\bigl(1+\log(2n)\bigr)\cdot k\log 2\cdot 6k^{2}\log k\cdot k\notag,\\
			&> -5.4\cdot 10^{12} k^{6} (\log k)^{2} \log n, \label{eq:matveev-bound}
		\end{align}
		where we have used the fact that $1+\log k \leq 3\log k$ and  \begin{align*}
			1+\log {(2n)}&=1+\log 2 + \log n\\&=\left(\frac{1}{\log n}+\frac{\log 2}{\log n}+1\right)\log n\\& \leq 3\log n,
		\end{align*} for $n\geq k+2, $ and $k\geq 2.$
		Comparing \eqref{eq:Gamma-bound} with \eqref{eq:matveev-bound}, we obtain
		\[
		n\log\alpha - \log 4 < 5.4\cdot 10^{12} k^{6} (\log k)^{2} \log n,
		\]
		and therefore
		\[
		n < 1.2\cdot 10^{13} k^{6} (\log k)^{2} \log n.
		\]
		Proceeding via Lemma~\ref{Guz}, fix $x:=n, r:= 1$ and
		$T := 1.2\cdot 10^{13} k^{6} (\log k)^{2}$. Since $T>x/\left(\log x\right)^r$ and $T>(4r^2)^r,$ we have that 
		\begin{align*}
			n &< 2\bigl(1.2\cdot 10^{13} k^{6} (\log k)^{2}\bigr)
			\log\bigl(1.2\cdot 10^{13} k^{6} (\log k)^{2}\bigr),\\
			&< 1.3\cdot 10^{15} k^{6} (\log k)^{3}.
		\end{align*}
		This completes the proof of Lemma \ref{lem:n-bound}.
	\end{proof}
	We proceed to examine the two cases basing on $k$.
	\subsubsection{The case $k\leq500$}
	Suppose that \(k\leq500\). By Lemma \ref{lem:n-bound}, we have
	$$
	n<1.3\cdot 10^{15}k^6(\log k)^3<4.9\cdot 10^{33}.
	$$
	We now refine the upper bound for \(n\). For this purpose, we return to \eqref{eq:Gamma-bound} and write
	\[
	\Gamma_1
	=
	2^{\mathfrak p}
	\left(3f_k(\alpha)\alpha^n\right)^{-1}-1
	=
	e^{\Lambda_1}-1.
	\]
	We have already established that \(\Gamma_1\neq0\), and hence \(\Lambda_1\neq0\). Assume that \(n\geq3\). Then
	$$
	|e^{\Lambda_1}-1|=|\Gamma_1|<0.5,
	$$
	and inequality \eqref{eq.2,7} gives $ |\log(1+\Gamma_1)|<1.5|\Gamma_1|.$
	Therefore, we obtain
	\begin{align*}
     \Lambda_1=\mathfrak{p}\log2-n\log\alpha-\log(3f_k(\alpha)),
	\end{align*} which can be rearranged to \begin{align*}
	\left|n\frac{\log \alpha}{\log 2}-\mathfrak{p}+\frac{\log{(3f_k(\alpha))}}{\log 2}\right|<\frac{6}{\log 2}\cdot \alpha^{-n}.
	\end{align*}
	Applying Lemma \ref{DP} with the quantities \begin{align*}
		\tau &:= \frac{\log \alpha}{\log 2}, & \mu &:= \frac{\log {(3 f_k(\alpha))}}{\log 2}, &C &:= \frac{6}{\log 2},  & \beta &:= \alpha.
	\end{align*}Here $u=n$ and $v=\mathfrak{p}$. We set $M:=9.8\cdot10^{33}$, which is an upper bound for $2n$.
	For each $k\in[2,500]$ we computed $\alpha=\alpha(k)$, $\tau$ and $\mu$ to high precision
	and took the first convergent $r/s$ of $\tau$ with $s>6M$ and $\theta>0$. If a convergent
	gave $\theta\le 0$, we passed to the next one.  For $k=2$ this is the convergent $s_{68}$, with
	$s\approx 7.0157\cdot10^{34}$ and $\theta=0.00202508$, which gives
	\[
	n<\frac{\log\!\big((6/\log 2)\,s/\theta\big)}{\log\alpha}<100.52<101 .
	\]
	Over all $k\in[2,500]$ we found $\theta>4.2\cdot10^{-4}$, the smallest value occurring at $k=105$.
	Writing $H_k$ for the resulting bound on $n$, the interval $[k+2,\lfloor H_k\rfloor]$ is
	non-empty only for $k\le 86$. It therefore remains to check the equation $P_n^{(k)}=(2^{\mathfrak{p}}+1)/3$
	for $2\le k\le 86$, $n\in[k+2,101]$ and $\mathfrak{p}\in[3,2n]$. A search in SageMath~10.6 found no solutions. We therefore turn to the case $k>500$.
	
	\subsubsection{The case $k>500$}
	Suppose that \(k>500\). By Lemma \ref{lem:n-bound}, we have $$
	n<1.3\cdot 10^{15}k^6(\log k)^3<\varphi^{k/2},
	$$
	where the last inequality holds for all \(k>500\). Thus, $	n<\varphi^{k/2}, $ which allows us to apply inequality \eqref{pk_b1}. Combining this with \eqref{eq2.6m}, we obtain \begin{align*}
		\left|\frac{2^{\mathfrak{p}}+1}{3}-\frac{\varphi^{2n}}{\varphi+2}\right|&=\left|P_n^{(k)}-\frac{\varphi^{2n}}{\varphi+2}\right|,\\ 
		&\leq \left|P_n^{(k)}-f_k(\alpha)\alpha^n\right|+\left|f_k(\alpha)\alpha^n-\frac{\varphi^{2n}}{\varphi+2}\right|,\\
		&<\frac{1}{2}+\left(\frac{\varphi^{2n}}{\varphi +2}\right)\cdot \frac{4}{\varphi^{\frac{k}{2}}}.
	\end{align*}
	Therefore \begin{align*}
		\left|\frac{2^{\mathfrak{p}}}{3}-\frac{\varphi^{2n}}{\varphi+2}\right|&<\frac{1}{2}+\left(\frac{\varphi^{2n}}{\varphi +2}\right)\cdot \frac{4}{\varphi^{\frac{k}{2}}}+\frac{1}{3},\\
		&=\frac{5}{6}+\left(\frac{\varphi^{2n}}{\varphi +2}\right)\cdot \frac{4}{\varphi^{\frac{k}{2}}}.
	\end{align*}
	Dividing through by $\left(\frac{\varphi^{2n}}{\varphi +2}\right)$ yields \begin{align*}
		\left|\frac{2^\mathfrak{p}}{3}\cdot \left(\frac{\varphi+2}{\varphi^{2n}}\right)-1\right|<\frac{5}{6}\cdot \left(\frac{\varphi+2}{\varphi^{2n}}\right)+\frac{4}{\varphi^{\frac{k}{2}}}<\frac{4}{\varphi^{2n}}+\frac{4}{\varphi^{k/2}}.
	\end{align*}
	Since $2n>n\geq k+2>k/2$, then $\min\{\varphi^{2n},\varphi^{k/2}\}=\varphi^{k/2}$ and thus we get \begin{align*}
		\left|\frac{2^\mathfrak{p}}{3}\cdot \left(\frac{\varphi+2}{\varphi^{2n}}\right)-1\right|<\frac{8}{\varphi^{k/2}}.
	\end{align*}
	Fix \begin{align}\label{gam2}
		\left|\Gamma_2\right|:=\left|\frac{2^\mathfrak{p}}{3}\cdot \left(\frac{\varphi+2}{\varphi^{2n}}\right)-1\right|<\frac{8}{\varphi^{k/2}}.
	\end{align} Clearly $\Gamma_2\neq 0$ otherwise we would have \[2^\mathfrak{p}\left(\frac{\varphi+2}{3}\right)=\varphi^{2n},\] taking Norms $N_{\mathbb{K}/\mathbb{Q}}$ on both sides for $\mathbb{K}=\mathbb{Q}(\sqrt{5})$, we have \[N(2^\mathfrak{p}\left(\frac{\varphi+2}{3}\right))=N(\varphi^{2n}).\] Using the fact that $N\left(\varphi+2\right)=5$ and $N\left(\varphi^{2n}\right)=\left(\varphi\varphi'\right)^{2n}=1,$ we have that \[\frac{2^{2\mathfrak{p}}\cdot 5}{9}=1,\] which is impossible for $\mathfrak{p}\geq 3.$ \\Taking the algebraic Number field $\mathbb{K}=\mathbb{Q}(\sqrt{5})$, we have $D:=2, t=3$ and  \begin{align*}
	\gamma_1 &:= 2, & \gamma_2 &:= \varphi, & \gamma_3 &:= \frac{\varphi +2}{3},\\
	b_1 &:= \mathfrak{p}, & b_2 &:= -2n, & b_3 &:= 1.
	\end{align*} Moreover, we have
	$$
	h(\gamma_1)=h(2)=\log 2,\qquad
	h(\gamma_2)=h(\varphi)=\frac{1}{2}\log\varphi,
	$$
	and
	$$
	h(\gamma_3)
	=h\left(\frac{\varphi+2}{3}\right)
	<\frac{1}{2}\log\varphi+\log 12.
	$$
	Therefore, we take
	$$
	A_1:=2\log 2,\qquad
	A_2:=\log\varphi,\qquad
	A_3:=\log\varphi+2\log 12.
	$$ Also, since $B\geq \max\{|b_i|:i=1,2,3\}$, we take $B:=2n$ as before. Proceeding via Theorem \ref{matveev} gives \begin{equation}\label{gamma2}
		\begin{aligned}
			\log{\left|\Gamma_2\right|}
			&>-1.4\cdot 30^6\cdot 3^{4.5}\cdot 2^2(1+\log 2)(1+\log {2n})
			(2\log2)(\log \varphi)(\log \varphi+2\log {12})\\
			&>-7.1\cdot10^{12}\log n.
		\end{aligned} 
	\end{equation} Comparing \eqref{gam2} to \eqref{gamma2}, we get \begin{align*}
	\frac{k}{2}\log \varphi -\log 8 <7.1\cdot 10^{12}\log n,
	\end{align*} thus we obtain the bound \begin{align}\label{v2}
	k<3.03\times 10^{13}\log n.
	\end{align} Recall the bound for $n$ in inequality \eqref{eq:n-bound}. Substituting this bound in \eqref{v2} yields a new bound for $k$ as \begin{align}
	k<3.03\cdot 10^{13}\log\left(1.3\cdot 10^{15}k^6(\log k)^3\right)<4.5\cdot 10^{14}\log k.
	\end{align} Proceeding via Lemma \ref{Guz}, we fix $x:=k, r:=1$ and $T:=4.5\cdot 10^{14}.$ We therefore have that $T>x/\left(\log x\right)^r,$ and $T>(4r^2)^r$ which yields \begin{align}
	k<2(4.5\cdot 10^{14})(\log(4.5\cdot 10^{14}))<3.1\cdot10^{16}.
	\end{align} Applying this new bound for $k$ into inequality \eqref{eq:n-bound} gives a new bound for $n$ as $n<6.4\cdot10^{118}$.
	We proceed to reduce these bounds. To do this, we revisit \eqref{gam2} and rewrite it as \begin{align*}
	\left|\Gamma_2\right|:=\left|2^\mathfrak{p}\cdot \varphi^{-2n}\cdot \left(\frac{\varphi+2}{3}\right)-1\right|<\frac{8}{\varphi^{k/2}}.
	\end{align*} We already showed that $\Gamma_2\neq 0,$ it follows that $\Lambda_2\neq 0.$ Assuming $k\geq 12,$ then $\left|\Gamma_2\right|<8\varphi^{-k/2}<0.5,$ we obtain the inequality $\left|e^{\Lambda_2}-1\right|=\left|\Gamma_2\right|<0.5$ and by inequality \eqref{eq.2,7} we get $  \left|\log{(1+\Gamma_2)}\right|<1.5\left|\Gamma_2\right|$. Consequently, we derive \begin{align*}
	\left|\mathfrak{p}\log 2-n\log {\varphi^2}+\log {\left(\frac{\varphi+2}{3}\right)}\right|<\frac{12}{\varphi^{k/2}}
	\end{align*} which when rearranged yields \begin{align*}
	\left|\mathfrak{p}\frac{\log 2}{\log{\varphi^2}}-n+\frac{\log {\left(\frac{\varphi+2}{3}\right)}}{{\log {\varphi^2}}}\right|<\frac{12}{\log {\varphi^2}}\cdot \varphi^{-k/2}
	\end{align*}
	Applying Lemma \ref{DP} with the quantities \begin{align*}
		\tau &:= \frac{\log 2}{\log \varphi^2}, & \mu &:= \frac{\log {\left(\frac{\varphi+2}{3}\right)}}{\log \varphi^2}, &C &:= \frac{12}{\log {\varphi^2}},  & \beta &:= \sqrt{\varphi}.
	\end{align*}
	Let $\tau:=[a_0;a_1, a_2, \cdots ]=[0; 1, 2, 1, 1, 2, 1, 6, 1, 8, \cdots]$ be the continued fraction expansion of $\tau.$  Here $u=\mathfrak{p}$, $v=n$ and $\ell=k$. We set $M:=1.28\cdot10^{119}$, which is an upper bound for $2n$. The first convergent of $\tau$ with $s>6M$ and $\theta>0$ is $s_{234}$, with
	$s\approx1.22\cdot10^{120}$ and $\theta\approx0.3202$. Hence  \begin{align}\label{pk}
		k<\frac{2\log{\left((12/\log {\varphi^2}) s/\theta\right)}}{\log \varphi}<1165.
	\end{align} We repeat this reduction procedure but with a new value of $M:=2.3\cdot 10^{36}$ evaluated for the reduced bound $k$ in \eqref{pk} and find that $s=s_{78}>6M$ with $\theta=0.23119.$ This gives $k<379$ which contradicts our original assumption that $k>500$ and thus the proof is complete.
	\section*{Acknowledgments}
	The author gratefully acknowledges Makerere University for providing a conducive environment for this research. This work is dedicated to the memory of Professor Florian Luca (1969--2026), whose remarkable contributions to Number Theory profoundly inspired the author's research in Diophantine equations and recurrence sequences. It is offered as a modest tribute to his enduring mathematical legacy.
	\section*{Disclosure statement}
	The author reports no conflict of interest.
	
	\section*{Addresses}
	$ ^{1} $ Department of Mathematics, Makerere University, Kampala, Uganda.
	
	Email: \url{kaggwaprosper58@gmail.com}
	
	\end{document}